\documentclass[preprint]{elsarticle}
\usepackage{mathrsfs}
\usepackage{amsmath,amssymb}
\usepackage{lineno,hyperref}
\modulolinenumbers[5]

\newtheorem{thm}{Theorem}
 \newtheorem{lem}[thm]{Lemma}
  \newtheorem{cor}[thm]{Corollary}
    \newtheorem{conj}[thm]{Conjecture}
 \newtheorem{defn}[thm]{Definition}
 \newtheorem{prop}[thm]{Proposition}
 \newdefinition{rmk}{Remark}
 \newproof{pf}{Proof}
 \newproof{poti}{Proof of Theorem \ref{thm1}}
 \newproof{potii}{Proof of Theorem \ref{thm2}}

\begin{document}

\begin{frontmatter}

\title{The complete splittings of finite abelian groups}


\author[mymainaddress]{Kevin Zhao}
\ead{zhkw-hebei@163.com}

\address[mymainaddress]{Department of Mathematics, South China normal university, Guangzhou 510631, China}

\begin{abstract}

  Let $G$ be a finite group.
  We will say that $M$ and $S$ form a \textsl{complete splitting} (\textsl{splitting}) of $G$ if every element (nonzero element) $g$ of $G$ has a unique representation of the form $g=ms$ with $m\in M$ and $s\in S$, and $0$ has a such representation (while $0$ has no such representation).

  In this paper, we determine the structures of complete splittings of finite abelian groups.
  In particular, for complete splittings of cyclic groups our description is more specific.
  Furthermore, we show some results for existence and nonexistence of complete splittings of cyclic groups
  and find a relationship between complete splittings and splittings for finite groups.
\end{abstract}


\end{frontmatter}

\section{Introduction}
The splittings of finite abelian groups are closely related to lattice tilings
and are easy to be generated to lattice packings.

Let $K_0$ be a polytope composed of unit cubes and $v+K_0$ be a translate of $K_0$ for some vector $v$.
A family of translations $\{v+K_0: v\in H\}$ is called an integer lattice packing if $H$ is an $n$-dimensional subgroup of $\mathbb{Z}^n$ and, for any two vectors $v$ and $w$ in $H$, the interiors of $v + K_0$ and $w + K_0$ are disjoint;
furthermore, if the $n$-dimensional Euclidean space $\mathbb{Z}^n$ is contained in the family of these translations,
then we call it an integer lattice tiling.

The splitting problem can be traceable to a geometric problem posed by H. Minkowski \cite{[HM]} and solved by
G. Haj\'{o}s \cite{[GH]}.
This problem is closely related to the \textit{factorizations} of finite abelian groups introduced by G. Haj\'{o}s \cite{[GH]}.
Let $G$ be a finite abelian group, written additively, and let $A_1,\ldots,A_n$ be nonempty subsets of $G$.
If for each $g$ in $G$ there are unique elements $a_1,\ldots,a_n$ of $G$ such that
$$g= a_1 +\ldots+ a_n, \ a_1\in A_1,\ldots,a_n\in A_n,$$
then we say that $G= A_1+\ldots+A_n$ is a \textit{factorization} of $G$.
One can find more results about it in \cite{[SS]}.
If $G$ is written multiplicatively, then we also call  $G= A_1\cdot\ldots\cdot A_n$ a \textit{factorization} of $G$.

Stein \cite{St67} first studied the splitting problem and showed its equivalence to the problem of tiling the Euclidean space by translates of certain polytope composed of unit cubes.
Whereafter Stein and Hickerson etc. continued to study the splittings of finite groups.
More results can be found in \cite{GS81, HS74, [H], SS94, Sz86, Sz87}.
This problem attracted recent attention again due to their equivalence to codes correcting single limited magnitude errors in flash memories (see \cite{BE13, CSBB10, EB10, KBE11, KLNY11, KLY12, M96,  OSW18, Sc12, Sc14, YKB13, ZG16, ZG18, ZZG17} and the references therein).
For existence and nonexistence results on the splittings of finite abelian groups, one can refer to   \cite{KLNY11, KLY12, Sc12, Sc14, St67, W95, YKB13},  and \cite{[PK], ZG16, ZG18, ZZG17}.
Furthermore, some people studied the packings of finite abelian groups and the relationship between the lattice packing and the coding theory. For these problems one can refer to \cite{HS86, St84, YKB13}, and \cite{M96, St84, SS94}.
However, the determination of splitting (packing) structures of finite abelian groups is a wide open question in general.
Motivated by lattice tilings and lattice packings, we consider these problems in a new way.

\begin{defn}
Let $G$ be a finite group, written additively, $M$ a set of integers, and $S,\ K$  subsets of $G$.
We will say that $M$ and $S$ form a \textsl{partial splitting} of $G$ for $K$ if every element $g$ of $G\setminus K$ has a unique representation of the form $g=ms$ with $m\in M$ and $s\in S$.
Denote it by $G\setminus K=MS$.
If $0\in K$, then we can also call it a splitting of $G\setminus K$.

If $K=\{0\}$, we call it a \textsl{splitting} of $G$.

If $\{0\}\subseteq K$, we call it a \textsl{packing} of $G$.

If $K=\phi$, we call it a \textsl{complete splitting} of $G$.

If $K\neq\phi$, we call it a \textsl{proper partial splitting} of $G$.

$M$ will be referred to as the multiplier set and $S$ as the partial splitting (splitting, packing, complete splitting, proper partial splitting, respectively) set.
We will also say that $M$ partially splits (splits, packs, completely splits, partially splits, respectively) $G$ with partial splitting (splitting, packing, complete splitting, proper partial splitting, respectively) set $S$, or simply that $M$ partially splits (splits, packs, completely splits, partially splits, respectively) $G$, if the particular set $S$ is not of interest.
\end{defn}

We are interested in $K=\phi$ or $pG$ with prime $p$.
For a cyclic group $\mathbb{Z}_n$ of order $n$, a splitting of $\mathbb{Z}_n$ can imply a partial splitting of $\mathbb{Z}_n$ for $K=p\mathbb{Z}_n$.
It follows since
$$\{g\in \mathbb{Z}_n: gcd(g,p)=1\}=\{m\in M: gcd(m,p)=1\}\cdot \{s\in S: gcd(s,p)=1\}$$
if $G\setminus \{0\}=MS$ is a splitting.

By imitating the proof of Theorem 1.2.1 in \cite{[H]}, one can obtain the following theorem:

\begin{thm}\label{thmsubsp}
Let $H$ be a normal subgroup of a finite group $G$. ($H$ and $G$ need not be abelian.)
Suppose that $M$ completely splits $G$ and $M$ splits $G/H$.
Then every complete splitting $G=MS$ induces a complete splitting of $H$, i.e., $H= M(S\cap H)$.
\end{thm}

By Theorem \ref{thmsubsp}, it is easy to see that the following corollary holds:

\begin{cor} \label{sp-comsp}
Let $H$ be a subgroup of a finite abelian group $G$.
Suppose that $M$ completely splits $G$ with $G=MS$ and $0=mg$ with $m\in M$ and $g\in S$.
If $g\not \in H$,
then $M$ does not split $G/H$.
\end{cor}

\pf
Suppose $M$ splits $G/H$, set $G/H\setminus \{0\}=MT$.
Since $M$ completely splits $G$, by Theorem \ref{thmsubsp}, we have the compete splitting of $H$: $$H=M(S\cap H).$$
Since $0\in H$ and $m\in M$, we have $g\in S\cap H$, a contradiction.
\qed

\textbf{Remark:} The value of Corollary \ref{sp-comsp} is little.
Suppose that $M$ splits $G/H$.
Then we must have that $|M|$ divides $|G/H|-1$.
It follows that $gcd(|M|,|G/H|)=1$.
Since $M$ completely splits $G$, one can obtain that $|M|$ divides $|G|=|H|\cdot |G/H|$.
Thus $|M|$ divides $|H|$.
From Theorem \ref{abelian-thm} it follows that $|M|=ord(g)$.
Thus $g$ must be contained in $H$ and Corollary \ref{sp-comsp} is true.
However, if one obtains a similar result for a proper partial splitting of $G$ instead of the complete splitting,
we will think that it is of great significance.
In particular, for the case $K=pG$.
Unfortunately, it fails.
The main cause of failure is that we can not obtain a similar result as Theorem \ref{thmsubsp} for a proper partial splitting of $G$ for $K=pG$.

\section{Preliminaries}

In this section we provide notations used throughout this work, and main conclusions obtained. The following notations are fixed throughout this paper.

Let $(G,+,0)$ be a finite abelian group and $G=MS$ be a complete splitting.
Thus $|G|=1$ or $G = \mathbb{Z}_{n_1} \oplus \mathbb{Z}_{n_2}\oplus \ldots
\oplus \mathbb{Z}_{n_k}$ with $1 < n_1|n_2| \ldots |n_k$,
where $r(G)=k$ is the rank of $G$ and the exponent
$\exp(G)$ of $G$ is $n_k$.
Without loss of generation, assume that $M\subseteq \mathbb{Z}_{n_k}$ is an integer subset modular $n_k$.
If $|M|=1$ or $|S|=1$, then we call the complete splitting \textit{trivial}.
For a nontrivial complete splitting, it is easy to see that $0\not\in M$ and $0\not\in S$.
It follows that if $0=ms$ with $m\in M$ and $s\in S$, then $m\neq 0$ and $g\neq 0$.
For an integer subset $A$, denote by $gcd(A)$ the greatest common divisor of all elements of $A$ and
denote by $lcm(A)$ the least common multiple of all elements of $A$.
In particular, if $A=\{m,\ n\}$, we can denote $gcd(m,n)$ by $(m,n)$ for short.
For positive integers $n$ and $g$ with $(n,g)=1$,
let $ord_n(g)$ denote the minimal positive integer $l$ such that $g^l\equiv 1$ (mod $n$).

For any $m\in \mathbb{Z}$, denote by $\mathbb{Z}_m$ an additive cyclic group of order $m$,
denote by $C_m$ a multiplicative cyclic group of order $m$,
and let $\mathbb{Z}_m^*=\{g\in \mathbb{Z}_m:(g,m)=1\}$.

This paper mainly studies the complete splittings of finite abelian groups.
In particular, for cyclic groups we determine some cases of existence and nonexistence of their complete splittings.
Furthermore, we think that our results can be generated to finite nonabelian groups, and propose the following conjecture:
\begin{conj}
Let $G$ be an additive finite group (not need be abelian) with $0\in G$. If $G=MS$ is a complete splitting and $0=mg$ with $m\in M$ and $g\in S$, then $mS=mG$ and $Mg=<g>$, i.e., $M$ is a complete set of representatives modulo $ord(g)$.
\end{conj}

Our main results are the followings:

\begin{thm}\label{abelian-thm}
Let $G= \mathbb{Z}_{n_1} \oplus \mathbb{Z}_{n_2}\oplus \ldots \oplus \mathbb{Z}_{n_k}$ with $1 < n_1|n_2 \ldots |n_k$ be a finite abelian group and $0\in G$. If $G=MS$ is a complete splitting and $0=mg$ with $m\in M$ and $g\in S$, then $mS=mG$, $g=(0,0,\ldots ,y_k\frac{n_k}{(m,n_k)})$ with $gcd(y_k,(m,n_k))=1$ and $Mg=<g>$ with $|M|=(m,n_k)$, i.e., $$S=\{(g_j,i+y_{ij}\frac{n_k}{(m,n_k)}):g_j\in \mathbb{Z}_{n_1}\oplus \ldots \oplus \mathbb{Z}_{n_{k-1}},i\in [0,\frac{n_k}{(m,n_k)}-1] \ and \ y_{ij}\in \mathbb{Z}\}$$ and $M$ is a complete set of representatives modulo $ord(g)$.
\end{thm}

By Theorem \ref{abelian-thm}, it is easy to show that the complete splittings of finite abelian $p$-groups are trivial.

\begin{prop}\label{p^a-complete-splitting-exist}
Let $\alpha$ be a positive integer, $p$ be a prime and
let $(G,+,0)$ be a finite abelian $p$-group of order $p^{\alpha}$.
Suppose that there exists a complete splitting $G=MS$.
Then $|M|=1$, if $r(G)\geq 2$; $|M|=1$ or $|S|=1$, if $r(G)=1$.
\end{prop}

\pf
Suppose $|M|>1$, $|S|>1$ and let $G = \mathbb{Z}_{p^{\alpha_1}} \oplus \mathbb{Z}_{p^{\alpha_2}}\oplus \ldots
\oplus \mathbb{Z}_{p^{\alpha_k}}$ with $1 \leq \alpha_1\leq \alpha_2\leq \ldots \leq \alpha_k$ and $|G|=p^{\alpha}$.
From Theorem \ref{abelian-thm} it follows that
there exist a positive integer $\beta$ and two nonzero elements $m\in M$, $g\in S$ such that $0=mg$ satisfying that
$$g=(0,0,\ldots ,y_k\frac{p^{\alpha_k}}{(m,p^{\alpha_k})})\neq 0\ with \ (y_k,(m,p^{\alpha_k}))=1,$$ $$1<|M|=(m,p^{\alpha_k})=p^{\beta}<p^{\alpha_k},$$
$$mS=mG\ and \ Mg=<g>.$$
Thus
$$m=m_1p^{\beta}\ and \ g=(0,0,\ldots ,y_kp^{\alpha_k-\beta})$$
with $(m_1,p)=1$, $(y_k,p)=1$ and $1\leq\beta\leq \alpha_k-1$.
For $mS=mG$ we have that $$m\cdot(S\setminus \{g\})=mS\setminus \{0\}=mG\setminus \{0\}=p^{\beta}G\setminus \{0\}.$$
Since $(0,0,\ldots ,y_kp^{\alpha_k-1})=p^{\beta-1}g=p^{\beta}\cdot (0,0,\ldots ,y_kp^{\alpha_k-1-\beta})\neq 0$,
from $Mg=<g>$ it follows that
$$(0,0,\ldots ,y_kp^{\alpha_k-1})\in (<g>\setminus\{0\})\cap (p^{\beta}G\setminus \{0\})=(M\setminus \{m\})\cdot g\cap m\cdot (S\setminus \{g\}).$$
This is in contradiction to $(M\setminus \{m\})\cdot g\cap m\cdot (S\setminus \{g\})=\phi$.
Hence, either $|M|=1$, $S=G$ or $M=G$, $|S|=1$.
For the latter case, the rank $r(G)$ of $G$ must be $1$, i.e., if $r(G)\geq 2$, then $|M|=1$.
We complete the proof.

\qed

By Theorem \ref{abelian-thm}, we also obtain a more concise conclusion on the complete splittings of cyclic groups:

\begin{cor}\label{cyclic-cor}
Let $G$ be a cyclic group with order $n\in \mathbb{N}$.
If $G=MS$ is a complete splitting and $0=mg$ with $m\in M$ and $g\in S$, then
$M$ is a complete set of representatives modulo $\frac{n}{(g,n)}=(m,n)$ and
$S$ is a complete set of representatives modulo $\frac{n}{(m,n)}=(n,g).$
\end{cor}

\pf
By Theorem \ref{abelian-thm}, we have that
$$M=\{i+k_i(m,n):i\in [0,(m,n)-1] \ and \ k_i\in \mathbb{Z}\}$$
and
$$S=\{j+y_j\frac{n}{(m,n)}:j\in [0,\frac{n}{(m,n)}-1] \ and \ y_j\in \mathbb{Z}\}.$$
Thus $M$ is a complete set of representatives modulo $(m,n)$ and
$S$ is a complete set of representatives modulo $\frac{n}{(m,n)}$.
Similarly, by exchanging $M$ and $S$ in the above,
one can obtain that $S$ is a complete set of representatives modulo $(g,n)$ and
$M$ is a complete set of representatives modulo $\frac{n}{(g,n)}$.
This proof is complete.
\qed

\section{The proof of Theorem \ref{abelian-thm}}

\textit{The proof of Theorem \ref{abelian-thm}:}
Suppose $G=\mathbb{Z}_{n_1}\oplus \mathbb{Z}_{n_2}\oplus \ldots \oplus \mathbb{Z}_{n_k}$ with $n_1|n_2|\ldots|n_k$.
Let $\varphi:$ $G\rightarrow mG$ be a homomorphism with $ker(\varphi)=\{g_0\in G:mg_0=0\}$ and $\varphi(g_0)=mg_0$ for any $g_0\in G$.
Since $G=MS$ is a complete splitting and $0=mg$ with $m\in M$ and $g\in S$, we have that $|G|=|M|\cdot |S|$, $mS\subseteq mG$ and $Mg\subseteq \{g_0\in G:mg_0=0\}=ker(\varphi)=\{(y_1\frac{n_1}{(m,n_1)},y_2\frac{n_2}{(m,n_2)},\ldots ,y_k\frac{n_k}{(m,n_k)}):
y_i\in \mathbb{Z} \ for \ 1\leq i\leq k\}$.
It follows that

$$
|mG|=\frac{|G|}{(m,n_1)(m,n_2)\ldots (m,n_k)}\geq |mS|=|S|,
$$
and

$$
|ker(\varphi)|=|\{g_0\in G:mg_0=0\}|=(m,n_1)(m,n_2)\ldots (m,n_k)\geq |Mg|=|M|.
$$
Therefore, $|G|=|M|\cdot |S|\leq |mG|\cdot |ker(\varphi)|=|G|$, which means that

$$
mS=mG, |S|=\frac{|G|}{(m,n_1)(m,n_2)\ldots (m,n_k)},
$$
and
\begin{gather}
Mg=\{g_0\in G:mg_0=0\},|M|=(m,n_1)(m,n_2)\ldots (m,n_k).   \label{eq1}
\end{gather}
From (\ref{eq1}) it is easy to see that $g=(y_1\frac{n_1}{(m,n_1)},y_2\frac{n_2}{(m,n_2)},\ldots , y_k\frac{n_k}{(m,n_k)})$
with $(y_i,(m,n_i))=1$ for $1\leq i\leq k$.

For any $g_0=(x_1,x_2, \ldots ,x_k)\in G$, denote by $(g_0)_i$ the $i$-th coordinate $x_i$ of $g_0$.
Set $G_{0i}=\{g_0\in G:mg_0=0, (g_0)_i=0\}$, $M_i=\{m_0\in M:(m_0g)_i=0\}$
and $d_i=gcd\{m_0:m_0\in M_i\}$.
It is easy to see that $M_i=\{m_0\in M:(n_i,m)|m_0\}$, $(n_i,m)|d_i$ and
\begin{gather}
d_iG \supseteq M_iS.\label{eq2}
\end{gather}
In addition, for $\{g_0\in G:mg_0=0\}=Mg$ we must have that $G_{0i}=M_ig$.
Thus
\begin{gather}
|M_i|=|G_{0i}|=\frac{(m,n_1)(m,n_2)\ldots (m,n_k)}{(m,n_i)}.\label{eq3}
\end{gather}
Combining (\ref{eq2}) and (\ref{eq3}) yields that
\begin{align*}
&|d_iG|=\frac{|G|}{(d_i,n_1)(d_i,n_2)\ldots (d_i,n_k)}\geq|M_i|\cdot |S|\\
&=\frac{(m,n_1)(m,n_2)\ldots (m,n_k)}{(m,n_i)}\cdot \frac{|G|}{(m,n_1)(m,n_2)\ldots (m,n_k)}\\
&=\frac{|G|}{(m,n_i)}.\\
\end{align*}
It follows that $$(m,n_i)\geq (d_i,n_1)(d_i,n_2)\ldots (d_i,n_k).$$
Putting $i=k$, we have that $(n_k,m)|d_k$, and
$(m,n_k)\geq (d_k,n_1)(d_k,n_2)\ldots (d_k,n_k)\geq (d_k,n_1)(d_k,n_2)\ldots ((m,n_k),n_k)\geq (m,n_k)$.
It follows that $(d_k,n_k)=(m,n_k)$ and $(d_k,n_j)=1$ for $1\leq j\leq k-1$.
For $(n_k,m)|d_k$ and $(n_j,m)|(n_k,m)$ for $1\leq j\leq k$, we have $(n_j,m)|(d_k,n_j)$  for $1\leq j\leq k-1$.
Thus $(n_j,m)=1$ for $1\leq j\leq k-1$.
Hence, $|M|=(n_k,m)$, $|S|=\frac{|G|}{(n_k,m)}$ and

$$
g=(0,0,\ldots ,y_k\frac{n_k}{(m,n_k)})\in S
$$
with $(y_k,(m,n_k))=1$.
It follows that $|M|=ord(g)$ and
$$Mg=\{g_0\in G:mg_0=0\}=\{(0,0,\ldots ,a\cdot \frac{n_k}{(m,n_k)})\in G:a\in [0,(m,n_k)-1]\}=<g>.$$
Thus
$$M=\{m,1+k_1(m,n_k),\ldots ,((m,n_k)-1)+k_{(m,n_k)-1}(m,n_k):k_i\in \mathbb{Z}\ for\ 1\leq i\leq (m,n_k)-1\}$$
is a complete set of representatives modulo $ord(g)$.
In addition, it is easy to see that
$$mS=mG=\{(g_0,(m,n_k)g_1):g_0\in \mathbb{Z}_{n_1}\oplus \ldots \oplus \mathbb{Z}_{n_{k-1}},g_1\in [0,\frac{n_k}{(m,n_k)}-1]\}.$$
Thus
$$S=\{(g_j,i+y_{ij}\frac{n_k}{(m,n_k)}):g_j\in \mathbb{Z}_{n_1}\oplus \ldots \oplus \mathbb{Z}_{n_{k-1}},i\in [0,\frac{n_k}{(m,n_k)}-1] \ and \ y_{ij}\in \mathbb{Z}\}.$$
\qed

\section{The Propositions of Complete Splitting of finite Cyclic Group}

In Corollary \ref{cyclic-cor}, we have shown a result on the complete splitting $\mathbb{Z}_n=MS$ of cyclic group $\mathbb{Z}_n$.
In this section, we will continue to study the problem for determining the structures of $M$ and $S$.
In the following, we show some results for existence and nonexistence of the complete splitting.

\begin{prop}
Let $(G,+,0)$ be a finite cyclic group of order $n\in \mathbb{N}$ and $M$ be a subset of $\mathbb{N}$.
If $|M|=k$, $\{k-3,k-2,k-1,k\}\subseteq M$ and $lcm(k-3,k-2,k-1,k)|n$,
then $M$ does not completely split $G$.
\end{prop}
\textbf{Remark:} If, in the above proposition, we replace "$\{k-3,k-2,k-1,k\}\subseteq M$"
by "$\{m_0,m_1,m_2,m_3\}\subseteq M$ with $(m_i,n)=k-i$ for $0\leq i\leq 3$", the conclusion still holds.

\pf
Suppose $M$ completely splits $G$, then there exists a subset $S$ of $G$ such that $G=MS$ and every element $g$ of $G$ has a unique representation $g=ms$ with $m\in M$ and $s\in S$.
Since $|G|=n$, $|M|=k$, $k\in M$ and $lcm(k-3,k-2,k-1,k)|n$, we have that $|S|=\frac{n}{k}$ and $k\cdot S\subseteq <k>\subseteq G$.
From $|kS|=|S|=\frac{n}{k}=|<k>|$ it follows that $$kS=<k>.$$
Thus, $$(M\setminus \{k\})\cdot S=MS\setminus kS=MS\setminus <k>=G\setminus <k>.$$
For $k-1\in M$ we have that $(k-1)S\cup (<k-1>\cap <k>)\subseteq <k-1>$ and
$|(k-1)S|+|<k-1>\cap <k>|=|S|+|<lcm(k-1,k)>|=\frac{n}{k}+\frac{n}{k(k-1)}=\frac{n}{k-1}=|<k-1>|$.
Thus $$(k-1)S\cup (<k-1>\cap <k>)= <k-1>.$$
It follows that $$(M\setminus \{k,k-1\})\cdot S=G\setminus \{<k>,<k-1>\}.$$
For $k-2\in M$ we have that $(k-2)S\cup (<k-2>\cap \{<k>,<k-1>\})\subseteq <k-2>$.
Thus $|(k-2)S|+|<k-2>\cap \{<k>,<k-1>\}|\leq \frac{n}{k-2}$.
Furthermore,
\begin{align*}
&|(k-2)S|+|<k-2>\cap \{<k>,<k-1>\}|\\
&=|S|+|<k-2>\cap <k>|+|<k-2>\cap <k-1>|-|<k-2>\cap <k>\cap<k-1>|\\
&=\frac{n}{k}+\frac{n}{lcm(k,k-2)}+\frac{n}{lcm(k-1,k-2)}-\frac{n}{lcm(k,k-1,k-2)}\\
&=\frac{n}{k}+\frac{n(k,k-2)}{k(k-2)}+\frac{n}{(k-1)(k-2)}-\frac{n(k,k-2)}{k(k-1)(k-2)}\\
&=\frac{n}{k}+\frac{n}{(k-1)(k-2)}+\frac{n(k,k-2)}{k(k-1)}\\
&\geq \frac{n}{k}+\frac{n}{(k-1)(k-2)}+\frac{n}{k(k-1)}\\
&=\frac{n}{k-2}.\\
\end{align*}
Therefore, $$(k-2)S\cup (<k-2>\cap \{<k>,<k-1>\})= <k-2>,$$
and $(k,k-2)=1$, i.e., $2\nmid k.$
It follows that $$(M\setminus \{k,k-1,k-2\})\cdot S=G\setminus \{<k>,<k-1>,<k-2>\}.$$
For $k-3\in M$ we have that $(k-3)S\cup (<k-3>\cap \{<k>,<k-1>,<k-2>\})\subseteq <k-3>$.
Thus $|(k-3)S|+|<k-3>\cap \{<k>,<k-1>,<k-2>\}|\leq \frac{n}{k-3}$.
Furthermore,
\begin{align*}
&|(k-3)S|+|<k-3>\cap \{<k>,<k-1>,<k-2>\}| \\
&=|S|+|<k-3>\cap <k>|+|<k-3>\cap <k-1>|+|<k-3>\cap <k-2>| \\
&-|<k-3>\cap <k>\cap <k-1>|-|<k-3>\cap <k>\cap <k-2>|\\
&-|<k-3>\cap <k-1>\cap <k-2>|+|<k-3>\cap <k>\cap <k-1>\cap <k-2>| \\
&=\frac{n}{k}+\frac{n}{lcm(k,k-3)}+\frac{n}{lcm(k-1,k-3)}+\frac{n}{lcm(k-2,k-3)}-\frac{n}{lcm(k,k-1,k-3)} \\
&-\frac{n}{lcm(k,k-2,k-3)}-\frac{n}{lcm(k-1,k-2,k-3)}+\frac{n}{lcm(k,k-1,k-2,k-3)} \\
&=\frac{n}{k}+\frac{n(k,k-3)}{k(k-3)}+\frac{n(k-1,k-3)}{(k-1)(k-3)}+\frac{n}{(k-2)(k-3)} \\
&-\frac{n(k,k-3)(k-1,k-3)}{k(k-1)(k-3)}-\frac{n(k,k-2)(k,k-3)}{k(k-2)(k-3)}\\
&-\frac{n(k-1,k-3)}{(k-1)(k-2)(k-3)}+\frac{n(k,k-2)(k,k-3)(k-1,k-3)}{k(k-1)(k-2)(k-3)} \\
&=(k,k-3)(\frac{n}{k(k-3)}-\frac{2n}{k(k-1)(k-3)}-\frac{n}{k(k-2)(k-3)}+\frac{2n}{k(k-1)(k-2)(k-3)}) \\
&+\frac{n}{k}+\frac{2n}{(k-1)(k-3)}+\frac{n}{(k-2)(k-3)}-\frac{2n}{(k-1)(k-2)(k-3)} \\
&=\frac{n(k-3)(k,k-3)}{k(k-1)(k-2)}+\frac{n}{k}+\frac{n}{(k-2)(k-3)}+\frac{2n}{(k-1)(k-2)} \\
&> \frac{n}{k}+\frac{n}{k(k-1)}+\frac{n}{(k-1)(k-2)}+\frac{n}{(k-2)(k-3)}=\frac{n}{k-3}. \\
\end{align*}
This is a contradiction and $M$ does not completely split $G$.

\qed

Now we will study the existence of the complete splitting of cyclic group $\mathbb{Z}_n$.
We have shown that the complete splittings of finite abelian $p$-groups are trivial.
In the following, we find some nontrivial complete splittings for $n=pq$ where $p,\ q$ are distinct primes.
For $n=p^{\alpha}q^{\beta}$ one can obtain similar results by imitating the proof of the case $n=pq$.

\begin{lem} \label{order}
Let $m$, $n$, $d_1$, $d_2$ be positive integers, $(m,n)=1$.
If $g\in \mathbb{Z}_{mn}^*$ satisfies that $ord_{mn}(g)=\frac{\varphi(mn)}{d}$ with $d=(\varphi(m),\varphi(n))$,
$ord_{m}(g)=\frac{\varphi(m)}{d_1}$ and $ord_{n}(g)=\frac{\varphi(n)}{d_2}$,
then $d_1d_2|d$ and $(d_1,d_2)=1$.
\end{lem}

\pf
For $ord_{mn}(g)=\frac{\varphi(mn)}{d}$ we have that $g^{\frac{\varphi(mn)}{d}}\equiv 1$ (mod $mn$).
It follows that $g^{\frac{\varphi(mn)}{d}}\equiv 1$ (mod $m$) and $g^{\frac{\varphi(mn)}{d}}\equiv 1$ (mod $n$).
Thus $ord_{m}(g)|\frac{\varphi(mn)}{d}$ and $ord_{n}(g)|\frac{\varphi(mn)}{d}$, i.e.,
$lcm(ord_{m}(g),ord_{n}(g))|\frac{\varphi(mn)}{d}=ord_{mn}(g)$.
Set $a:=ord_{m}(g)$ and $b:=ord_{n}(g)$.
Thus $g^{a}\equiv 1$ (mod $m$) and $g^{b}\equiv 1$ (mod $n$).
It follows that $m|g^{ab}-1$ and $n|g^{ab}-1$.
For $(m,n)=1$ we have that $g^{ab}\equiv 1$ (mod $mn$) and then $ord_{mn}(g)|ab$.
Therefore,
$$lcm(ord_{m}(g),ord_{n}(g))=\frac{\varphi(mn)}{d}=\frac{\varphi(m)\cdot \varphi(n)}{d}.$$
Since $ord_{m}(g)=\frac{\varphi(m)}{d_1}$ and $ord_{n}(g)=\frac{\varphi(n)}{d_2}$,
we have that $\frac{\varphi(m)\cdot \varphi(n)}{d_1d_2}=ord_{m}(g)\cdot ord_{n}(g)=lcm(ord_{m}(g),ord_{n}(g))\cdot (ord_{m}(g),ord_{n}(g))=\frac{\varphi(m)\cdot \varphi(n)}{d}(\frac{\varphi(m)}{d_1},\frac{\varphi(n)}{d_2})$.
It follows that $$d=d_1d_2(\frac{\varphi(m)}{d_1},\frac{\varphi(n)}{d_2})=d_1d_2(\frac{\varphi(m)}{d}\cdot \frac{d}{d_1},\frac{\varphi(n)}{d}\cdot \frac{d}{d_2}).$$
For $d=(\varphi(m),\varphi(n))$, one can obtain that
$(\frac{d}{d_1},\frac{d}{d_2})|(\frac{\varphi(m)}{d}\cdot \frac{d}{d_1},\frac{\varphi(n)}{d}\cdot \frac{d}{d_2}).$
Set $(\frac{d}{d_1},\frac{d}{d_2})=x$ and
$(\frac{\varphi(m)}{d}\cdot \frac{d}{d_1},\frac{\varphi(n)}{d}\cdot \frac{d}{d_2})=x\cdot y.$
It follows that $d=d_1d_2xy=d_1d_2(\frac{d}{d_1},\frac{d}{d_2})\cdot y=d_1d_2(d_2xy,d_1xy)\cdot y=d_1d_2(d_1,d_2)xy^2$.
Thus $d_1d_2|d$ and $(d_1,d_2)y=1$, i.e., $(d_1,d_2)=1$ and $y=1$.

\qed

\begin{lem}[\cite{[SS]}, Theorem 7.1] \label{prime-factor}
Let $m$ and $n$ be relatively prime positive integers.
If $A=\{a_1,\ldots ,a_m\}$ and $B=\{b_1,\ldots ,b_n\}$ are sets of integers such that their sum set
$$A+B=\{a_i+b_j:1\leq i\leq m, 1\leq j\leq n\}$$
is a complete set of representatives $modulo$ $mn$,
then $A$ is a complete set of residues $modulo$ $m$ and $B$ is a complete set of residues $modulo$ $n$.
\end{lem}

\begin{lem}[\cite{[G]}] \label{integer-factor}
Let $n$ be a positive integer.
If $n=p_1^{\alpha_1}\cdot\ldots\cdot p_k^{\alpha_k}$,
then $$\mathbb{Z}_n^*\cong \mathbb{Z}_{p_1^{\alpha_1}}^*\times \ldots \times \mathbb{Z}_{p_k^{\alpha_k}}^*.$$
\end{lem}

\begin{lem}[\cite{[GGS]}, Lemma 3.2] \label{integer-factor-rank=2}
Let $n_1$, $n_2$ be positive integers.
Then $$C_{n_1}\times C_{n_1}=C_{(n_1,n_2)}\times C_{lcm(n_1,n_2)}.$$
\end{lem}

Let $\mathbb{Z}_{pq}^*=\{g\in \mathbb{Z}_{pq}:(g,pq)=1\}$ with distinct primes $p$ and $q$.
By Lemma \ref{integer-factor} and Lemma \ref{integer-factor-rank=2},
we have $$\mathbb{Z}_{pq}^*=\mathbb{Z}_{p}^*\times \mathbb{Z}_{q}^*=C_d\times C_{\frac{(p-1)(q-1)}{d}}=<x>\times <g>$$
where $d=(p-1,q-1)$ and $x$, $g$ are the generators of $C_d$, $C_{\frac{(p-1)(q-1)}{d}}$, respectively.
Set $$ord_{p}(g)=\frac{p-1}{d_1}\ and \ ord_{q}(g)=\frac{q-1}{d_2}.$$
Combining $(p,q)=1$ with Lemma \ref{order} yields that $d_1d_2|d$ and $(d_1,d_2)=1$.
Thus $d$ has a factorization $d=d'd''$ with $(d',d'')=1$ and $d_1|d'$, $d_2|d''$.
Continuing the above analysis we obtain the following lemma:

\begin{lem} \label{<g>-factor}
Suppose that $\mathbb{Z}_{\frac{(p-1)(q-1)}{d}}=A+B$
is a factorization with
$|A|=\frac{p-1}{d'}$ and $|B|=\frac{q-1}{d''}$. Set $M_2=\{g^a: a\in A\}$ and $S_2=\{g^b: a\in B\}$.
If $(d,\frac{p-1}{d})=1$ and $(d,\frac{q-1}{d})=1$,
then $$<g>=M_2S_2$$ is also a factorization
satisfying that $|M_2|=\frac{p-1}{d'}$, $|S_2|=\frac{q-1}{d''}$, $(|M_2|,|S_2|)=1$ and
all elements in $M_2$, $S_2$ are distinct modulo $p$, $q$, respectively.
\end{lem}

\pf
Since $d=(p-1,q-1)$, $(d,\frac{p-1}{d})=1$ and $(d,\frac{q-1}{d})=1$,
we have that
$$(\frac{p-1}{d'},\frac{q-1}{d''})=(\frac{p-1}{d}\times \frac{d}{d'},\frac{q-1}{d}\times \frac{d}{d''})=(d',d'')=1.$$
It follows that $$\mathbb{Z}_{\frac{(p-1)(q-1)}{d}}=A+B$$ is a factorization with
$|A|=\frac{p-1}{d'}$, $|B|=\frac{q-1}{d''}$ and $(|A|,|B|)=1$.
By Lemma \ref{prime-factor} we have that
$A$ is a complete set of residues $modulo$ $\frac{p-1}{d'}$
and $B$ is a complete set of residues $modulo$ $\frac{q-1}{d''}$.
Since $M_2=\{g^a: a\in A\}$, $S_2=\{g^b: a\in B\}$ and
$g$ is a generator of $C_{\frac{(p-1)(q-1)}{d}}$, we have that $$<g>=M_2S_2$$ is a factorization with $|M_2|=|A|=\frac{p-1}{d'}$, $|S_2|=|B|=\frac{q-1}{d''}$ and $(|M_2|,|S_2|)=1$.
Suppose that there exist two distinct elements $a_1$, $a_2$ in $A$ such that
$g^{a_1}\equiv g^{a_2}$ (mod $p$).
Thus $g^{a_1-a_2}\equiv 1$ (mod $p$) and $a_1\not\equiv a_2$ (mod $\frac{p-1}{d'}$).
From $ord_{p}(g)=\frac{p-1}{d_1}$ it follows that $\frac{p-1}{d_1}|a_1-a_2$.
For $d_1|d'$ we have that $\frac{p-1}{d'}|\frac{p-1}{d_1}$, and then $\frac{p-1}{d'}|a_1-a_2$.
This is in contradiction to $a_1\not\equiv a_2$ (mod $\frac{p-1}{d'}$).
Hence, all elements in $M_2$ are distinct modulo $p$.
Similarly, we can show than all elements in $S_2$ are distinct modulo $q$ and this complete the proof.

\qed

Continue the analysis of Lemma \ref{<g>-factor}
and we obtain a proposition on the existence of the complete splitting of cyclic group $\mathbb{Z}_{pq}$.

\begin{prop} \label{pq-complete-splitting-exist}
Let $(G,+,0)$ be a finite cyclic group with order $n=pq$.
Suppose that $(d,\frac{p-1}{d})=1$ and $(d,\frac{q-1}{d})=1$.
If one of the following conditions holds:
\begin{description}
  \item[(1)] $ord_{p}(x)=d'$ or $ord_{q}(x)=d''$,
  \item[(2)] $d_1=1$ or $d_2=1$,
  \item[(3)] $d=p_0^{\alpha}$ with a prime $p_0$ and a positive integer $\alpha$,
\end{description}
then there exists integral subsets $M$ and $S$ such that $G=MS$ is a complete splitting.
\end{prop}

\pf
$(1)$ If $ord_{p}(x)=d'$,
then let $A=<\frac{q-1}{d''}>$ and $B=[0,\frac{q-1}{d''}-1]$ be two subsets of $\mathbb{Z}_{\frac{(p-1)(q-1)}{d}}$.
It follows that $\mathbb{Z}_{\frac{(p-1)(q-1)}{d}}=A+B$
is a factorization with
$|A|=\frac{p-1}{d'}$ and $|B|=\frac{q-1}{d''}$.
By Lemma \ref{<g>-factor} one has that
$$<g>=M_2S_2$$ with $M_2=\{g^a: a\in A\}$ and $S_2=\{g^b: a\in B\}$ is a factorization
satisfying that $|M_2|=\frac{p-1}{d'}$, $|S_2|=\frac{q-1}{d''}$, $(|M_2|,|S_2|)=1$ and
all elements in $M_2$ ($S_2$) are distinct modulo $p$ (modulo $q$), respectively.

We claim that
$$\mathbb{Z}_{pq}^*=M_1S_1$$ with $M_1=\cup_{i=0}^{d'-1}x^iM_2$ and $S_1=\cup_{j=0}^{d''-1}x^{d'j}S_2$
is a factorization satisfying that
$M_1$, $S_1$ are reduced residue systems $modulo$ $p$, $q$, respectively.

For $\mathbb{Z}_{pq}^*=<x>\times <g>=\cup_{i=0}^{d-1}x^i<g>$,
it is easy to see that $\mathbb{Z}_{pq}^*=M_1S_1$.
Since $|M_1|=d'|M_2|=p-1$ and $|S_1|=d''|S_2|=q-1$,
we have that $|\mathbb{Z}_{pq}^*|=(p-1)(q-1)=|M_1|\cdot |S_1|$.
It follows that $\mathbb{Z}_{pq}^*=M_1S_1$ is a factorization.

Suppose there exist two distinct elements $x^{i_1}g^{a_1},\ x^{i_2}g^{a_2}\in M_1$ such that
$$x^{i_1}g^{a_1}\equiv x^{i_2}g^{a_2}\ (mod \ p)$$ where $i_1$, $i_2\in [0,d'-1]$ and $g^{a_1}$, $g^{a_2}\in M_2$ with $a_1,\ a_2 \in \{\ell\cdot\frac{q-1}{d''}: \ell\in [0,\ \frac{p-1}{d'}-1]\}$.
Since all elements in $M_2$ are distinct modulo $p$,
we must have $i_1\neq i_2$.
It follows that $a_1\neq a_2$, since otherwise $x^{i_1-i_2}\equiv 1$ (mod $p$), and
this is in contradiction in $0<|i_1-i_2|<d'=ord_{p}(x)$.
Without loss of generation, assume that $i_1<i_2$.
Thus  there exists an integer $k=(a_1-a_2)/\frac{q-1}{d''}\in [-(\frac{p-1}{d'}-1),(\frac{p-1}{d'}-1)]^*$ such that
$$x^i\equiv g^{\frac{q-1}{d''}k} \ (mod \ p)$$
where $0<i=i_2-i_1<d'=ord_{p}(x)$.
For $ord_{p}(x)=d'$,
we have that $x^{id'}\equiv g^{\frac{q-1}{d''}kd'}\equiv 1$ (mod $p$).
From $ord_{p}(g)=\frac{p-1}{d_1}$ it follows that $\frac{p-1}{d_1}|\frac{q-1}{d''}kd'$,
that is $$\frac{p-1}{d}\times \frac{d}{d'}\times \frac{d'}{d_1}|\frac{q-1}{d}\times \frac{d}{d''}\times kd'.$$
Since $d=(p-1,q-1)$, $(d,\frac{p-1}{d})=1$, $(d,\frac{q-1}{d})=1$, $d=d'd''$ and $(d',d'')=1$,
we have that $(\frac{p-1}{d}\times \frac{d}{d'},\frac{q-1}{d}\times \frac{d}{d''}\times d')
=(\frac{p-1}{d}\times d'',\frac{q-1}{d}\times d'^2)=1$.
Therefore, $$\frac{p-1}{d}\times \frac{d}{d'}=\frac{p-1}{d'}|k.$$
This is in contradiction to $k\in [-(\frac{p-1}{d'}-1),(\frac{p-1}{d'}-1)]^*$.
Hence, all elements in $M_1$ are distinct modulo $p$.

Since $M_1\subseteq \mathbb{Z}_{pq}^*$, we have that $(m_1,p)=1$ for any $m_1\in M_1$.
For $|M_1|=p-1$ one can obtain that $M_1$ is a reduced residue system $modulo$ $p$.

Suppose there exist two distinct elements $x^{d'j_1}g^{b_1},\ x^{d'j_2}g^{b_2}\in S_1$ such that
\begin{gather}
x^{d'j_1}g^{b_1}\equiv x^{d'j_2}g^{b_2}\ (mod \ q)   \label{eqmod}
\end{gather}
where $j_1$, $j_2\in [0,d''-1]$ and $g^{b_1}$, $g^{b_2}\in S_2$ with $b_1,\ b_2 \in [0,\frac{q-1}{d''}-1]$.
Since all elements in $S_2$ are distinct modulo $q$,
we must have $j_1\neq j_2$. Without loss of generation, assume that $j_1<j_2$.
Set $t=ord_q(x)$.
For $ord_{p}(x)=d'$, we have that $x^{d'}\equiv 1$ (mod $p$).
Combining it with $x^{t}\equiv 1$ (mod $q$) yields that $x^{td'}\equiv 1$ (mod $pq$).
From $ord_{pq}(x)=d=d'd''$ it follows that $d|td'$, i.e., $d''|t$.
If $b_1=b_2$, then $x^{(j_1-j_2)d'}\equiv 1$ (mod $q$) with $|j_1-j_2|\in [1,d''-1]$.
It follows that $t|(j_1-j_2)d'$ and then $d''|(j_1-j_2)d'$.
For $(d',d'')=1$, we have that $d''|(j_1-j_2)$ and
this is in contradiction in $|j_1-j_2|\in [1,d''-1]$.
Therefore, $b_1\neq b_2$.
From (\ref{eqmod}) it follows that  $$x^{d'j}\equiv g^{b_0} \ (mod \ q)$$
where $j=j_2-j_1\in [1,d''-1]$ and $b_0=b_1-b_2\in [-(\frac{q-1}{d''}-1),(\frac{q-1}{d''}-1)]^*$.
For $x^d\equiv 1$ (mod $pq$),
we have that $x^{d'jd''}\equiv x^{dj}\equiv g^{b_0d''}\equiv 1$ (mod $q$).
From $ord_{q}(g)=\frac{q-1}{d_2}$ it follows that $\frac{q-1}{d_2}|b_0d''$,
that is $$\frac{q-1}{d}\times \frac{d}{d''}\times \frac{d''}{d_2}|b_0d''.$$
Since $(d,\frac{q-1}{d})=1$, $d=d'd''$ and $(d',d'')=1$,
it is easy to see that $(\frac{q-1}{d}\times \frac{d}{d''},d'')=1$.
Therefore, $$\frac{q-1}{d}\times \frac{d}{d''}=\frac{q-1}{d''}|b_0.$$
This is in contradiction to $b_0\in [-(\frac{q-1}{d''}-1),(\frac{q-1}{d''}-1)]^*$.
Hence, all elements in $S_1$ are distinct modulo $q$.

Since $S_1\subseteq \mathbb{Z}_{pq}^*$, we have that $(s_1,q)=1$ for any $s_1\in S_1$.
For $|S_1|=q-1$ one can obtain that $S_1$ is a reduced residue system $modulo$ $q$.
This complete the proof of the claim.

If $ord_{q}(x)=d''$,
then let $A=[0,\frac{p-1}{d'}-1]$, $B=<\frac{p-1}{d'}>$ be two subsets of $\mathbb{Z}_{\frac{(p-1)(q-1)}{d}}$.
By imitating the proof of the claim,
one can obtain that
$$\mathbb{Z}_{pq}^*=M_1S_1$$ with $M_1=\cup_{i=0}^{d'-1}x^{d''i}M_2$, $S_1=\cup_{j=0}^{d''-1}x^{j}S_2$
is a factorization satisfying that
$M_1$, $S_1$ are reduced residue systems $modulo$ $p$, $q$, respectively.

Let $$M=\{p\}\cup M_1 \ and\ S=\{q\}\cup S_1.$$
Thus it's easy to see that $G=MS$ is a complete splitting.

$(2)$ If $d_1=1$,
then let $d'=1$ and $d''=d$.
Take two subsets $A=<\frac{q-1}{d''}>=<\frac{q-1}{d}>$, $B=[0,\frac{q-1}{d''}-1]=[0,\frac{q-1}{d}-1]$ of $\mathbb{Z}_{\frac{(p-1)(q-1)}{d}}$
and repeat the reasoning of the case $ord_{p}(x)=d'$ in $(1)$.
One can show that $M_1=M_2$ is a reduced residue system $modulo$ $p$ and
$S_1=\cup_{j=1}^{d}x^{j}S_2$ is a reduced residue system $modulo$ $q$.

If $d_2=1$,
then let $d'=d$ and $d''=1$.
Take two subsets $A=[0,\frac{p-1}{d'}-1]=[0,\frac{p-1}{d}-1]$, $B=<\frac{p-1}{d'}>=<\frac{p-1}{d}>$ of $\mathbb{Z}_{\frac{(p-1)(q-1)}{d}}$
and repeat the reasoning of the case $ord_{q}(x)=d''$ in $(1)$.
One can show that $M_1=\cup_{i=1}^{d}x^{i}M_2$ is a reduced residue system $modulo$ $p$ and
$S_1=S_2$ is a reduced residue system $modulo$ $q$.

Thus $G=MS$ with $M=\{p\}\cup M_1$ and $S=\{q\}\cup S_1$ is a complete splitting.

$(3)$ Since $d=p_0^{\alpha}$ with prime $p_0$ and positive integer $\alpha$,
from $d_1d_2|d$ and $(d_1,d_2)=1$ it follows that $d_1=1$ or $d_2=1$.
By $(2)$, we complete the proof.

\qed

\textbf{Example:} Put $p=31$ and $q=43$. Thus $d=(p-1,q-1)=6$.
Let $x=6$ and $g=45$.
By a calculation one can show that $ord_{pq}(x)=d=6$, $ord_{pq}(g)=\frac{(p-1)(q-1)}{d}=210$,
$ord_{p}(g)=15=\frac{p-1}{2}$, $ord_{q}(g)=14=\frac{q-1}{3}$, $ord_{q}(x)=3$ and $<6>\cap <45>=\{1\}$.
It follows that $\mathbb{Z}_{pq}^*=<6>\times <45>$, $d_1=2$ and $d_2=3$.
Thus take $d'=2$ and $d''=3$.
Let $A=[0,\frac{p-1}{d'}-1]=[0,14]$ and $B=<\frac{p-1}{d'}>=<15>$ be two subsets of $\mathbb{Z}_{\frac{(p-1)(q-1)}{d}}=\mathbb{Z}_{210}$.
Let $M_2=\{g^a: a\in A\}$ and $S_2=\{g^b: a\in B\}$.
For $ord_{q}(x)=3=d''$, $(d,\frac{p-1}{d})=(6,5)=1$ and $(d,\frac{q-1}{d})=(6,7)=1$, by the proof of Proposition \ref{pq-complete-splitting-exist} $(1)$
we have that $G=MS$ is a complete splitting,
where $M=\{p\}\cup M_1$ and $S=\{q\}\cup S_1$ with
$M_1=\cup_{i=0}^{d'-1}x^{d''i}M_2=M_2\cup x^{3}M_2$ and $S_1=\cup_{j=0}^{d''-1}x^{j}S_2=S_2\cup xS_2\cup x^{2}S_2$.

\section{The proof of Theorem \ref{thmsubsp}}

For a subset $A$ of a group $G$, let $c(A)$ be the number of nonzero elements of $A$; that is, $c(A) =|A|$ if $0\not\in A$ and $c(A)=|A|-1$ if $0\in A$.

In \cite{[H]}, the author defined the $M$-partition of a finite group:
\begin{defn}[\cite{[H]}]
\begin{description}
  \item[(a)] A \textit{partition} of a set $X$ is a set of disjoint nonempty subsets of $X$ whose union is $X$.
If $\mathscr {G}$ is a partition of $X$, then the equivalence relation associated with $\mathscr {G}$
is denoted by "$\sim _{\mathscr {G}}$"; that is, $x\sim _{\mathscr {G}} y$ if $x\in A$ and $y\in A$ for some $A\in \mathscr {G}$. (Definition 1.0.0)
  \item[(b)] A partition $\mathscr {G}$ of a group $G$ is a $M$-\textit{compatible} provided that, for any $g$ and $h$ in $G$ and $m\in M$, if $g\sim _{\mathscr {G}} h$ then $mg\sim _{\mathscr {G}} mh$.
Equivalently, if $A\in \mathscr {G}$ and $m\in M$, then $mA\subseteq B$ for some $B\in \mathscr {G}$.

Given $A$, $B\in \mathscr {G}$, we let $q(A,B)=\{m\in M| mB\subseteq A\}$. (Definition 1.0.1)
  \item[(c)] Let $A$ and $B$ be elements of an $M$-compatible partition $\mathscr {G}$ of a group $G$. Then $B$ \textit{divides} $A$ (written "$B|A$") if there are elements $B = B_0, B_1, \ldots ,B_r =A$ of $\mathscr {G}$ such that, for $0\leq i<r$, $m_iB_i\subseteq B_{i+1}$ for some $m_i\in M$.
Equivalently, $B$ divides $A$ if $nB\subseteq A$ for some $n$ which can be expressed as a (possibly empty) product of elements of $M$.

We will say that $B$ is a \textit{proper divisor} of $A$ (written "$B<A$") if $B|A$ but $B\neq A$. (Definition 1.0.4)
  \item[(d)] An $M$-compatible partition $\mathscr {G}$ of a group $G$ is called an $M$-\textit{partition} of $G$ if divisibility is a partial ordering on $\mathscr {G}$. (Definition 1.0.5)
\end{description}
\end{defn}

In the proof of the main theorem, we use a special $M$-partition, i.e., the order partition:
\begin{defn}[\cite{[H]}, Definition 1.1.1]
Let $G$ be a finite group. The \textit{order partition} $\Theta$ of $G$ is defined by $g\sim _{\Theta}h$ if $o(g) = o(h)$. The \textit{order} of an equivalence class $A$ in $\Theta$ is the common order of the elements of $A$.
\end{defn}

To prove Theorem \ref{thmsubsp}, we need some lemmas and their generations to the complete splittings.

\begin{lem}[\cite{[H]}, Lemma 1.2.2] \label{lemH3}
Let $H$ be a normal subgroup of a finite group $G$ and $\mathscr {G}$ an $M$-partition of $G/H$.
For $A\in \mathscr {G}$, let $A^*=\{g\in G|g+H\in A\}$. Let $\mathscr {G}^*=\{A^*|A\in \mathscr {G}\}$.
Then, for $m\in M$ and $A$, $B\in \mathscr {G}$, we have:
\begin{enumerate}
  \item $\mathscr {G}^*$ is an $M$-partition of $G$;
  \item $mA\subseteq B$ if and only if $mA^*\subseteq B^*$;
  \item $A|B$ if and only if $A^*|B^*$;
  \item $|A^*|=|H|\cdot |A|$.
\end{enumerate}
That is, $\mathscr {G}$ and $\mathscr {G}^*$ have the same structure with respect to scalar multiplication, and each element of $\mathscr {G}^*$ is $|H|$ times as large as the corresponding element of $\mathscr {G}$.
\end{lem}

\begin{lem}[\cite{[H]}, Theorem 1.0.3] \label{lemH1}
Suppose $G$ is a finite group with a splitting $G\setminus \{0\}=MS$, $\mathscr {G}$ is an $M$-compatible partition of $G$, and $A\in \mathscr {G}$. Then

$$
c(A)=\sum_{B\in \mathscr {G}}|q(A,B)|\cdot |S\cap B|.
$$
\end{lem}

Repeat the reasoning in the proof of Lemma \ref{lemH1} and one can obtain the following lemma:
\begin{lem}
Suppose $G$ is a finite group with a complete splitting $G=MS$, $\mathscr {G}$ is a $M$-compatible partition of $G$, and $A\in \mathscr {G}$. Then

$$
|A|=\sum_{B\in \mathscr {G}}|q(A,B)|\cdot |S\cap B|.
$$
\end{lem}
\pf
Since $G=MS$ is a complete splitting and $\mathscr {G}$ is a $M$-compatible partition of $G$,
we have that for any two distinct elements $m_1$, $m_2\in M$ and $B$, $B_1\in \mathscr {G}$,
$$
(M(S\cap B_1))\cap (M(S\cap B))=\phi \ and \ (m_1(S\cap B))\cap (m_2(S\cap B))=\phi.
$$
In addition, for any $A$, $B\in \mathscr {G}$, $m\in M$ and $s\in S\cap B$, we have that $ms\in A$ if and only if $mB\subseteq A$ if and only if $m\in q(A,B)$. Thus,
$$|M(S\cap B)\cap A|=\sum_{m\in M}|m(S\cap B)\cap A|=\sum_{m\in q(A,B)}|m(S\cap B)|=|q(A,B)|\cdot |(S\cap B)|.$$
For $G=MS=\cup_{B\in \mathscr {G}}B$ we have that $S=S\cap G=\cup_{B\in \mathscr {G}}(S\cap B)$.
Combining the above results yields that
\begin{align*}
|A|&=|A\cap G|=|A\cap MS|=|A\cap M(\cup_{B\in \mathscr {G}}(S\cap B))|=|\cup_{B\in \mathscr {G}}(M(S\cap B)\cap A)| \\
&=\sum_{B\in \mathscr {G}}|M(S\cap B)\cap A|=\sum_{B\in \mathscr {G}}|q(A,B)|\cdot |S\cap B| \\
\end{align*}
\qed

The proof of the following lemma is exactly the same as one of Lemma 1.0.6 in \cite{[H]}.
\begin{lem} \label{nonempty}
Let $\mathscr {G}$ is an $M$-partition of a group $G$. Then
\begin{enumerate}
  \item If, for each $g\in G$, $M$ contains an element $m$ relatively prime to $o(g)$,
  then $q(A,A)\neq \phi$ for every $A\in \mathscr {G}$.
  \item If, for some $g$, $M$ has no such element, then $M$ does not (completely) split $G$.
\end{enumerate}
\end{lem}

\begin{lem}[\cite{[H]}, Theorem 1.0.8] \label{thmin1}
Let $G$ be a finite group and $\mathscr {G}$ an $M$-partition of $G$.
If $G$ has a splitting $G\setminus \{0\}=MS$.
Then the values of $|S\cap A|$ for $A\in \mathscr {G}$ can be computed recursively by the formula

$$
|S\cap A|=\frac{1}{|q(A,A)|}(c(A)-\sum_{B\in \mathscr {G},B<A}|q(A,B)|\cdot |S\cap B|),
$$
where the sum is over all $B\in \mathscr {G}$ for which $B<A$.
\end{lem}

Combining Lemma \ref{nonempty} and the proof of Lemma \ref{thmin1} yields the following lemma:
\begin{lem} \label{thmin2}
Let $G$ be a finite group and $\mathscr {G}$ an $M$-partition of $G$.
If $G$ has a splitting a complete splitting $G=MS$.
Then the values of $|S\cap A|$ for $A\in \mathscr {G}$ can be computed recursively by the formula

$$
|S\cap A|=\frac{1}{|q(A,A)|}(|A|-\sum_{B\in \mathscr {G},B<A}|q(A,B)|\cdot |S\cap B|),
$$
where the sum is over all $B\in \mathscr {G}$ for which $B<A$.
\end{lem}

By imitating the proof of Theorem 1.2.1 in \cite{[H]}, we can obtain Theorem \ref{thmsubsp}:

\textit{The proof of Theorem \ref{thmsubsp}:}
Let $\mathscr {G}$ be an $M$-partition of $G/H$ in which $0$ is in a class by itself; for example, we may take $\mathscr {G}$ to be the order partition. Let $\mathscr {G}^*$ be as defined in Lemma \ref{lemH3}.
Let $G=MS$ and $G/H\setminus \{0\}=MT$.

We must show that $H=M(S\cap H)$. Clearly, $M(S\cap H)\subseteq H$. Also, the products $ms$ are all distant where $m\in M$ and $s\in S$,
so $|M(S\cap H)|=|M|\cdot |S\cap H|$. Thus it is sufficient to show that

$$
|S\cap H|=\frac{|H|}{|M|}.
$$

By Lemma \ref{thmin2}, we have
\begin{gather}
|S\cap A^*|=\frac{1}{|q(A^*,A^*)|}(|A^*|-\sum_{B\in \mathscr {G}^*,B<A^*}|q(A^*,B)|\cdot |S\cap B|)
\label{|SA^*|}
\end{gather}
for any $A\in \mathscr {G}$.

We now claim that
$$
|S\cap A^*|=|H|\cdot |T\cap A|
$$
where $A\in \mathscr {G}$ and $A\neq \{0\}$.

The proof of the claim is by induction on the number $n$ of proper divisors of $A^*$.
It is easy to see that $n=|\{B\in \mathscr {G}^*|B<A^*\}|=|\{B'\in \mathscr {G}|B'<A\}|$.
If $n=0$, then from (\ref{|SA^*|}) it follows that
\begin{gather}
|S\cap A^*|=\frac{|A^*|}{|q(A^*,A^*)|}=\frac{|A|\cdot |H|}{|q(A,A)|}.  \label{|0SA^*|}
\end{gather}
Since $G/H\setminus \{0\}=MT$ and $A\in \mathscr {G}\setminus \{0\}$,
for any $a\in A$ there exist $m\in M$ and $t\in T$ such that $mt=a\in A$.
Thus $t\in A$, since otherwise there exists some $B\in \mathscr {G}\setminus A$ with $t\in B$,
and it would follow that $mB\subseteq A$, i.e., $B<A$ from the definition of $M$-compatible,
which is in contradiction to $n=0$.
It follows that $t\in T\cap A$ and by the definition of $M$-compatible one has that $mA\subseteq A$,
i.e., $m\in q(A,A)$.
Thus $A\subseteq q(A,A)\cdot (T\cap A)$.
It is easy to see that $q(A,A)\cdot (T\cap A)\subseteq A$, and then $A=q(A,A)\cdot (T\cap A)$.
From (\ref{|0SA^*|}) it follows that $|S\cap A^*|=|H|\cdot |T\cap A|$.

Now suppose that claim is true for all proper divisors of $A^*$.
Since each $B\in \mathscr {G}^*$ has the form $C^*$ for some $C\in \mathscr {G}$ and
$B|A^*$ if and only if $C|A$, combining them with (\ref{|SA^*|}) yields that
\begin{gather}
|S\cap A^*|=\frac{1}{|q(A^*,A^*)|}(|A^*|-\sum_{C\in \mathscr {G},C<A}|q(A^*,C^*)|\cdot |S\cap C^*|).  \label{|CSA^*|}
\end{gather}
Part ($2$) of Lemma \ref{lemH3} implies that $q(A^*,C^*)=q(A,C)$.
Part ($4$) implies that $|A^*|=|A|\cdot |H|$.
By the inductive hypothesis we have $|S\cap C^*|=|H|\cdot |T\cap C|$.
Combining these results with (\ref{|CSA^*|}) yields that
\begin{equation}
\label{|1CSA^*|}
   \begin{aligned}
   |S\cap A^*|&=\frac{1}{|q(A,A)|}(|A|\cdot |H|-\sum_{C\in \mathscr {G},C<A}|q(A,C)|\cdot |H|\cdot |T\cap C|) \\
   &=|H|\frac{1}{|q(A,A)|}(|A|-\sum_{C\in \mathscr {G},C<A}|q(A,C)|\cdot |T\cap C|).  \\
   \end{aligned}
  \end{equation}
Since $G/H\setminus \{0\}=MT$ is a splitting, from Lemma \ref{thmin1} it follows that
\begin{gather}
|T\cap A|=\frac{1}{|q(A,A)|}(c(A)-\sum_{C\in \mathscr {G},C<A}|q(A,C)|\cdot |T\cap C|).  \label{|TA|}
\end{gather}
For $0\not \in A$, we have $c(A)=|A|$.
Thus by (\ref{|1CSA^*|}) and (\ref{|TA|}) one can obtain that
$|S\cap A^*|=|H|\cdot |T\cap A|$ and the proof of the claim is complete.

Now let $A=\{0\}\in \mathscr {G}$. Thus $A^*=H$, $q(A,A)=M$ and $|T\cap A|=0$.
In addition, for any $C\in \mathscr {G}\setminus A$ we have $C\neq \{0\}$ and
then from the claim it follows that $|S\cap C^*|=|H|\cdot |T\cap C|$.
Part ($2$) of Lemma \ref{lemH3} implies that $q(A^*,C^*)=q(A,C)$ and $q(A^*,A^*)=q(A,A)$.
Combining these results with (\ref{|SA^*|}) yields that
\begin{equation}
\label{|CSA^*0|}
   \begin{aligned}
   |S\cap A^*|=|S\cap H|&=\frac{1}{|q(A^*,A^*)|}(|A^*|-\sum_{C\in \mathscr {G},C<A}|q(A^*,C^*)|\cdot |S\cap C^*|)\\
&=\frac{1}{|q(A,A)|}(|H|-\sum_{C\in \mathscr {G},C<A}|q(A,C)|\cdot |H|\cdot |T\cap C|) \\
&=|H|\frac{1}{|M|}(1-\sum_{C\in \mathscr {G},C<A}|q(A,C)|\cdot |T\cap C|). \\
   \end{aligned}
\end{equation}
By Lemma \ref{lemH1} it is easy to see that
\begin{equation}
\label{|CSA^*1|}
   \begin{aligned}
   c(A)=0&=\sum_{C\in \mathscr {G}}|q(A,C)|\cdot |T\cap C|\\
&=|q(A,A)|\cdot |T\cap A|+\sum_{C\in \mathscr {G},C<A}|q(A,C)|\cdot |T\cap C|\\
&=\sum_{C\in \mathscr {G},C<A}|q(A,C)|\cdot |T\cap C|.\\
   \end{aligned}
\end{equation}
By (\ref{|CSA^*0|}) and (\ref{|CSA^*1|}) we have that

$$
|S\cap H|=\frac{|H|}{|M|}
$$
and the proof is complete.

\qed

\textbf{Open Problem:} Let $\mathbb{Z}_n$ be a finite cyclic group. Prove a similar result as Theorem \ref{thmsubsp} for a proper partial splitting of $G$ for $K=pG$.
\section*{Acknowledgments}
This work is supported by NSF of China  (Grant No. 11671153).
The author is sincerely grateful to professor Pingzhi Yuan for his guidance and the anonymous referee
for useful comments and suggestions.

\section*{References}


\end{document}